\documentclass[12pt]{article}
\usepackage{epsfig,amsmath,amsthm,amssymb,graphics,amsfonts,psfrag}
\usepackage[letterpaper,margin=1in]{geometry}

\def\R{\mathbb{R}}

\def\P{\mathbb{P}}
\def\E{\mathbb{E}}

\def\P{\mathbb{P}}

\newcommand{\be}{\begin{equation}}
\newcommand{\ee}{\end{equation}}
\newcommand{\bea}{\begin{eqnarray}}
\newcommand{\eea}{\end{eqnarray}}
\newcommand{\beann}{\begin{eqnarray*}}
\newcommand{\eeann}{\end{eqnarray*}}
\newcommand{\benn}{\begin{equation*}}
\newcommand{\eenn}{\end{equation*}}

\def\ra{\rightarrow}
\def\I{\infty}
\def\I{\infty}

\newcommand{\cA}{{\mathcal A}}  % calligraphic A
\newcommand{\cB}{{\mathcal B}}  % calligraphic B
\newcommand{\cE}{{\mathcal E}}  % calligraphic E
\newcommand{\cH}{{\mathcal H}}  % calligraphic H
\newcommand{\cN}{{\mathcal N}}  % calligraphic N
\newcommand{\cO}{{\mathcal O}}  % calligraphic O
\begin{document}
 
\title{Deterministic continuation of stochastic metastable equilibria via Lyapunov equations and ellipsoids}
\author{Christian Kuehn\thanks{Max Planck Institute for the Physics of Complex Systems}}

\maketitle

\begin{abstract}
Numerical continuation methods for deterministic dynamical systems have been one of the most successful tools in applied dynamical systems theory. Continuation techniques have been employed in all branches of the natural sciences as well as in engineering to analyze ordinary, partial and delay differential equations. Here we show that the deterministic continuation algorithm for equilibrium points can be extended to track information about metastable equilibrium points of stochastic differential equations (SDEs). We stress that we do not develop a new technical tool but that we combine results and methods from probability theory, dynamical systems, numerical analysis, optimization and control theory into an algorithm that augments classical equilibrium continuation methods. In particular, we use ellipsoids defining regions of high concentration of sample paths. It is shown that these ellipsoids and the distances between them can be efficiently calculated using iterative methods that take advantage of the numerical continuation framework. We apply our method to a bistable neural competition model and a classical predator-prey system. Furthermore, we show how global assumptions on the flow can be incorporated - if they are available - by relating numerical continuation, Kramers' formula and Rayleigh iteration.
\end{abstract}

{\bf Keywords:} Numerical continuation, bifurcation analysis, metastability, stochastic dynamics, covariance, Lyapunov equation, ellipsoids, iterative methods, neural competition, predator-prey system, Rayleigh iteration, Kramers' law. 

\section{Introduction}  
\label{sec:intro}

Consider a deterministic dynamical system given by a differential equation
\be
\label{eq:DE}
\frac{\partial x}{\partial t}=x'=L(x;\mu)
\ee
where $x$ represents phase space variables, $\mu\in\R$ is a parameter and $L$ is an operator or a map that describes a deterministic equation e.g.~an ordinary differential equation (ODE), partial differential equation (PDE) or delay differential equation (DDE). Time-independent solutions of \eqref{eq:DE} with $x'=0$ are steady states (or equilibria) $x^*=x^*(\mu)$ with $L(x^*(\mu);\mu)=0$. Given an equilibrium $x^*(\mu_1)$, numerical continuation allows us to efficiently compute how it changes under parameter variation i.e.~to compute $x^*(\mu_2)$ for small $|\mu_1-\mu_2|$. In the case of an ODE we have a vector field
\benn
L(x;\mu)=f(x;\mu)\qquad \text{with $f:\R^n\times \R\ra \R^n$.}
\eenn 
Numerical continuation can be used to compute a curve of equilibrium points $\gamma=\{x=x^*(\mu)\}$ which solves the algebraic equations $f(x;\mu)=0$. Furthermore, one can compute so-called test (or bifurcation) functions for each point on this curve that indicate a change of stability of the equilibrium point under parameter variations.

Introductions to numerical continuation can be found in \cite{Doedel,Kuznetsov,Govaerts,AllgowerGeorg}. There are also many software packages available with various standard continuation algorithms and test functions such as MatCont \cite{MatCont,DhoogeGovaertsKuznetsov}, AUTO \cite{Doedel_AUTO2007,KrauskopfOsingaGalan-Vioque}, PyDSTool \cite{ClewleySherwoodLaMarGuckenheimer} and DDE-BIFTOOL \cite{EngelborghsLuzyaninaSamaey}. The literature on the applications of numerical continuation techniques is extremely large. For example, it can be used to compute periodic and homoclinic orbits \cite{Kuznetsov}, stable and unstable invariant manifolds of equilibrium points \cite{KrauskopfOsinga1}, slow manifolds \cite{GuckenheimerKuehn2} and canard orbits \cite{DesrochesKrauskopfOsinga2} in fast-slow systems as well as isochrons \cite{OsingaMoehlis}, just to name a few. Application areas range from physics \cite{GreenKrauskopfSamaey,MercaderBatisteAlonsoKnobloch}, chemistry \cite{Koper,DoedelHeinemann} and biology \cite{KuznetsovFeoRinaldi,WechselbergerWeckesser} to engineering \cite{SeydelHlavaceka,RankinDesrochesKrauskopfLowenberg}. It is even possible to implement continuation methods directly in experiments \cite{Sieberetal}.  

Despite this success story, there seems to be very little work to extend continuation ideas to stochastic differential equations (SDEs). Current numerical approaches to SDEs mostly focus on simulation and forward integration \cite{KloedenPlaten,MilsteinTretyakov}. Other available methods are set-valued techniques \cite{DellnitzJunge} to track invariant measures and the direct solution of forward or backward Kolmogorov PDEs \cite{SpencerBergman,SchenkHoppe1}. An approach that tries to utilize classical continuation for stochastic problems is the moment map formulation \cite{BarkleyKevrekidisStuart,ErbanKevrekidisAdalsteinssonElston} where the primary motivation seems to arise from equation-free modelling \cite{MakeevMaroudasKevrekidis}.\\

However, suppose we have already used numerical continuation for a deterministic ODE and found stable equilibrium points or more general stable invariant sets. Then it is a natural question to ask how small noise influences the stability of these objects. In general, we expect a change to metastable invariant sets \cite{ArnoldSDE,BerglundGentz5} so that noise-induced transitions between different stable sets can occur. In this paper, we show that there is a very natural and straightforward extension of equilibrium continuation in the context of SDEs that provides local information about metastable equilibrium points. Our approach can be applied during a numerical continuation calculation or, slightly less efficiently, as a post-processing tool.\\ 

\textit{Remark:} We note that the algorithm we develop here is expected to extend to much wider classes of problems such as nonstationary solutions \cite{KuehnSDEcont2} as well as SPDEs \cite{PrevotRoeckner} and SDDEs \cite{ReissRiedleVanGaans}.\\

The method is based on combining well-known results and numerical techniques from different areas of mathematics and computing. A reader interested in getting an overview of our main steps should consider the analytical example presented in Section \ref{sec:ana_ex}. The general development based on minimal local assumptions is presented in Sections \ref{sec:meta_SDE}-\ref{sec:main}. We test our approach for a planar vector field modelling neuronal competition in Section \ref{sec:num_ex1}. In this example, we focus on the algorithmic performance and show how to integrate the algorithm in standard numerical continuation software. In Section \ref{sec:num_ex2} we consider the Rosenzweig-MacArthur predator-prey system and demonstrate that important dynamical systems conclusions and direct interpretations for applications can be obtained from our computational framework. Further examples of how our algorithm relates to important conclusions regarding the dynamics of a system can be found in \cite{BerglundGentz,BerglundGentzKuehn}. In Section \ref{sec:Kramer} a special case with a global gradient-structure assumption is considered.\\ 

\textit{Preliminary Remark 1:} All computations have been carried out in MatLab \cite{MatLab2010b}, version R2010b on a standard quad-core 2.4 GHz CPU with 4 GB RAM. The numerical continuation calculations of deterministic equilibrium points use version 2.5.1. of cl$\_$MatCont \cite{MatCont}.\\ 

\textit{Preliminary Remark 2:} All norms refer to the Euclidean norm so that we simply use the notation $\|\cdot\|$ instead of $\|\cdot\|_2$. All vectors are assumed to be column vectors. The superscript notation $(~)^T$ will denote the transpose of vectors/matrices and $I$ is going to denote an identity matrix of suitable size for the algebraic operation considered.

%%%%%%%%%%%%%%%%%%%%%%%%%%%%%%%%%%%%%%%%%%%%%%%%%%%%%%%%%%%%%%%%%%%%%%%%%%%%%%%%%%%%%%
\section{An Analytical Example}
\label{sec:ana_ex}

\begin{figure}[htbp]
\psfrag{u}{$u$}
\psfrag{x}{$x$}
\psfrag{xs}{$x^*$}
\psfrag{xp}{$x^+$}
\psfrag{xm}{$x^-$}
\psfrag{mu}{$\mu$}
\psfrag{sigma}{$\sigma$}
\psfrag{u1}{$U_{\text{\scriptsize{$-0.5$}}}$}
\psfrag{u2}{$U_{\text{\scriptsize{$0.5$}}}$}
\psfrag{u3}{$U_{\text{\scriptsize{$1.5$}}}$}
	\centering
		\includegraphics[width=0.95\textwidth]{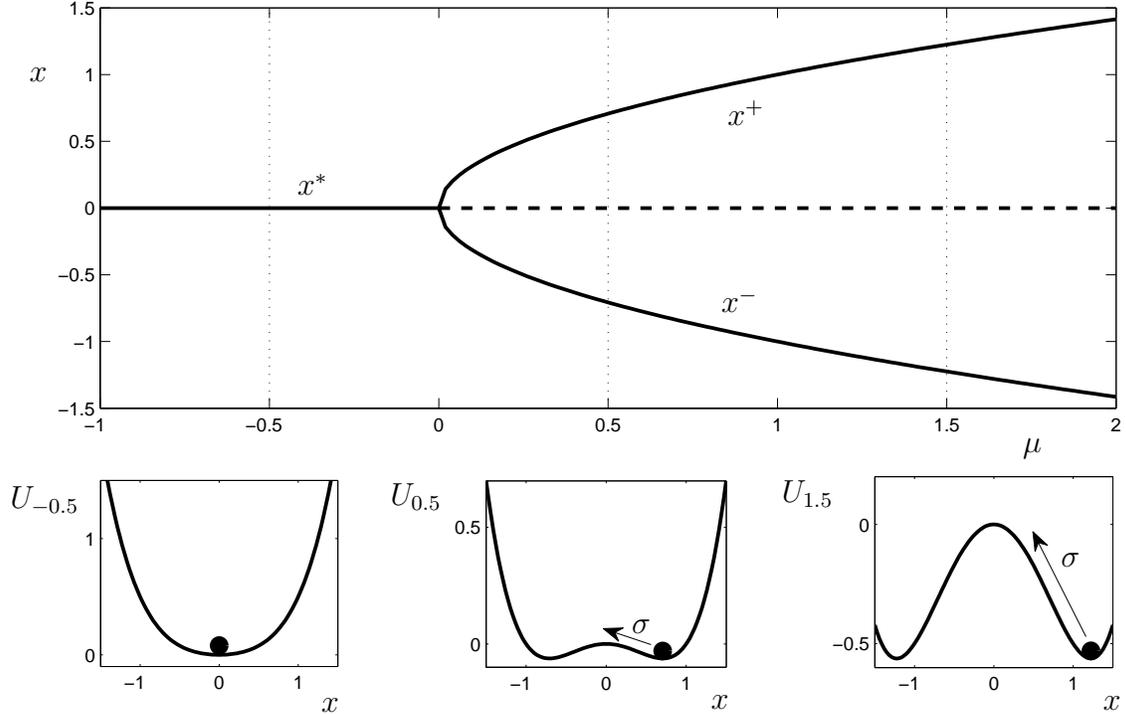}
	\caption{\label{fig:fig1}Bifurcation diagram (top) for the deterministic part of \eqref{eq:simple_ex}-\eqref{eq:simple_ex1}. The potentials are shown as well (bottom). In the bistable regime for $\mu>0$ and $\sigma>0$ noise-induced transitions between the metastable equilibria occur.}
\end{figure}

We start with a well-known analytical example to motivate the type of problems we are interested in and to present the basic conceptual ideas for the numerical analysis. Consider the following 1-dimensional SDE with additive noise 
\be
\label{eq:simple_ex}
dx_t=\left(\mu x_t-x_t^3\right)dt+\sigma dW_t=:f(x_t;\mu)dt+\sigma dW_t
\ee 
where $W_t$ is standard Brownian motion \cite{Oksendal}, $\sigma$ controls the noise level and $\mu\in\R$ is the main bifurcation parameter. Systems of the form \eqref{eq:simple_ex} appear very frequently in applications ranging from mean-field and Ising-type models for phase transitions in classical physics \cite{PikovskyZaikindelaCasa,BrokateSprekels}, reaction-rate theory in chemistry \cite{HaenggiTalknerBorkovec,PollakTalkner}, single neuron modelling \cite{Lindneretal} and bistable ecosystems \cite{GuttalJayaprakash2} in biology. The deterministic part of the SDE is a normal form for a pitchfork bifurcation \cite{Kuznetsov,Golubitsky2}. The dynamics of \eqref{eq:simple_ex} is easily understood by writing it as a gradient system
\be
\label{eq:simple_ex1}
dx_t=-\nabla U_\mu(x)dt+\sigma dW_t\qquad \text{with $U_\mu(x):=-\frac{\mu}{2}x^2+\frac{1}{4}x^4$.}
\ee 
so that the stochastic process $x_t$ can be interpreted as a particle moving in a potential $U_\mu(x)$. There is always one trivial deterministic equilibrium for \eqref{eq:simple_ex1} given by $x=x^*=0$. For $\mu<0$ the equilibrium $x^*$ is globally attracting and corresponds to a unique minimum of the potential $U_\mu(x)$. At $\mu=0$ a pitchfork bifurcation occurs; see Figure \ref{fig:fig1}. The equilibrium $x^*$ is destabilized and becomes a local maximum (saddle point) of the potential and two new locally stable equilibria $x^\pm=\pm\sqrt\mu$ appear for $\mu>0$ corresponding to minima of $U_\mu(x)$. Although we can easily obtain the deterministic equilibrium curves given in Figure \ref{fig:fig1} analytically as $\{(x,\mu)\in\R^2:x=0\}$ and $\{(x,\mu)\in\R^2:x=\pm\sqrt\mu\}$ one has to use numerical techniques, such as numerical continuation, for more general systems.

\begin{figure}[htbp]
\psfrag{a}{(a)}
\psfrag{b}{(b)}
\psfrag{c}{(c)}
\psfrag{d}{(d)}
\psfrag{t}{$t$}
\psfrag{x}{$x$}
\psfrag{mu}{$\mu$}
\psfrag{B+}{$\cB^+(h)$}
\psfrag{B-}{$\cB^-(h)$}
	\centering
		\includegraphics[width=1\textwidth]{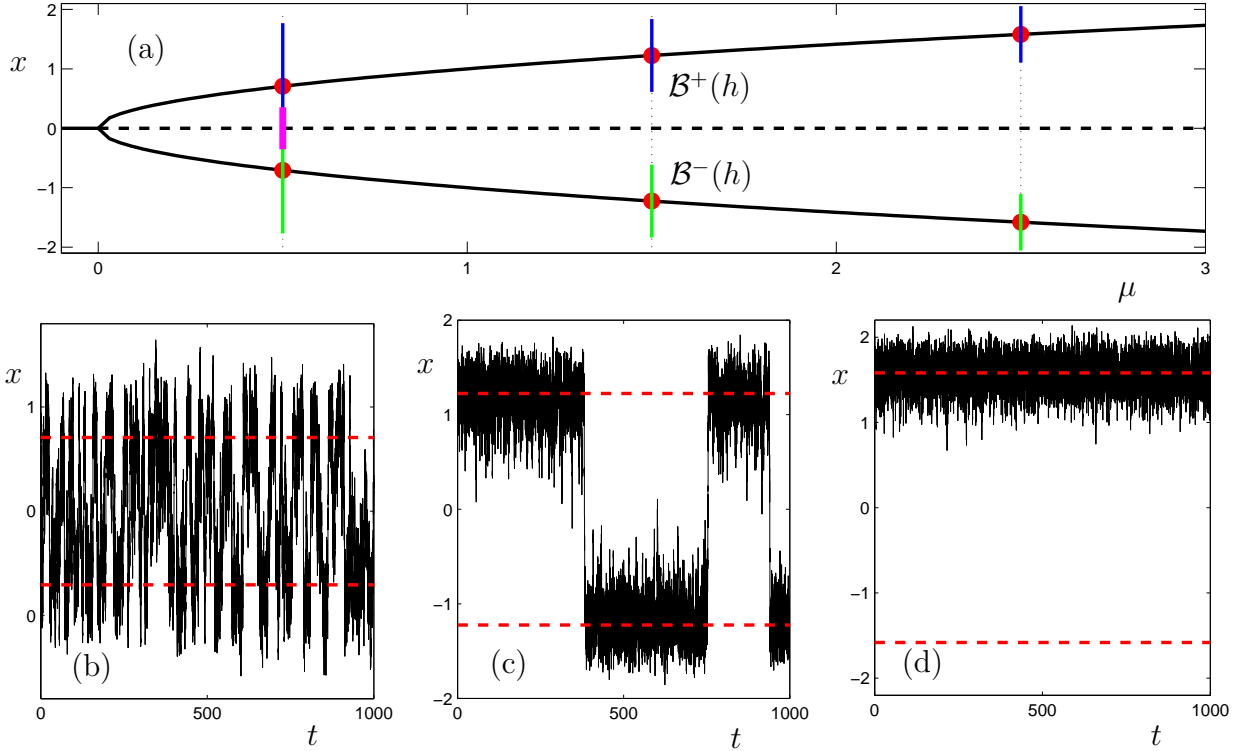}
	\caption{\label{fig:fig2}(a) Bifurcation diagram with stochastic neighbourhoods of the metastable equilibrium points. The deterministic equilibrium curves are shown in black. The stochastic variance neighbourhoods $\cB^\pm(h)$ defined by \eqref{eq:Bh} are indicated in blue/green with confidence level $h=6$. The overlap of $\cB^+(h)$ and $\cB^-(h)$ for $\mu=0.5$ is marked by a thick line (magenta). (b)-(d) Time series of \eqref{eq:simple_ex} with noise level $\sigma=0.5$ for (b) $\mu=0.5$, (c) $\mu=1.5$ and (d) $\mu=2.5$. The dashed curves (red) indicate the metastable equilibria $x^\pm$.}
\end{figure}

Interesting noise-induced dynamics occurs in the bistable regime for $\mu>0$. Fix any $\mu>0$, $\sigma>0$ and initial condition $x_0$. Then consider the first hitting times $t^\pm:=\inf\{t\geq 0:x_t=x^\pm\}$. A standard result from probability \cite{FreidlinWentzell} is that 
\be
\label{eq:meta_p1}
\P(t^\pm<\I)=1
\ee
i.e.~no matter where we start, we will eventually visit both deterministically stable equilibrium points with probability one. Although the result \eqref{eq:meta_p1} is of importance from a theoretical viewpoint it is of very limited practical use. In particular, the time scale on which the stochastic switching between the potential minima occurs is of major interest. Suppose we start the process $x_t$ at $x_0=x^+$. If $\sigma\gg 1$ frequent switching occurs and we will quickly visit $x^-$ while for $0<\sigma \ll 1$ switching is rare; see Figure \ref{fig:fig2}. The theory of large deviations \cite{FreidlinWentzell} considers the first-exit time over the saddle point $x^*$ given by $\tau^+:=\inf\{t\geq 0:x_0=x^+,x_t<x^*\}$ and shows that the mean first exit time is
\be
\label{eq:Kramers}
\E[\tau^+]=\cO\left(e^{2[U_\mu(x^*)-U_\mu(x^+)]/\sigma^2}\right) \qquad \text{as $\sigma\ra 0$.}
\ee
The result \eqref{eq:Kramers} is also known as Arrhenius' law \cite{Arrhenius} and the rate $1/\E[\tau^+]$ is called Eyring-Kramers rate \cite{Eyring,Kramers}; see also Section \ref{sec:Kramer}. Furthermore observe that the potential difference in \eqref{eq:Kramers} is given by
\benn 
U_\mu(x^*)-U_\mu(x^+)=0+\frac{\mu}{2}\sqrt\mu^2-\frac{1}{4}\sqrt\mu^4=\frac{\mu^2}{4}.
\eenn
Hence the switching probability/rate also depends on the bifurcation parameter $\mu$ and increases when $\mu\ra 0^+$. In Figure \ref{fig:fig2} we show three time series for a fixed noise level with varying bifurcation parameter $\mu>0$ over a fixed time interval $t\in[0,1000]$. It is clear that the dynamics in Figure \ref{fig:fig2}(b) with very frequent stochastic switching is different from rare switching events in Figure \ref{fig:fig2}(c) and no switching events up to $t=1000$ in Figure \ref{fig:fig2}(d).

One possibility to capture the stochastic behaviour is to solve the forward Kolmogorov (or Fokker-Planck) equation 
\cite{Oksendal} associated with \eqref{eq:simple_ex} given by
\be
\label{eq:FP}
\frac{\partial }{\partial t}p(x,t)=-\frac{\partial}{\partial x}(f(x;\mu)p(x,t))+\frac{\sigma^2}{2}\frac{\partial^2}{\partial x^2}p(x,t)
\ee   
where $p(x,t)=p(x,t|x_0,t_0)$ denotes the transition probability density of the stochastic process $x_t$ starting from $x_0$ at time $t_0$. However, solving \eqref{eq:FP} essentially solves the SDE \eqref{eq:simple_ex} everywhere in phase space. It is clear that for higher-dimensional nonlinear problems - where we are only interested in the local metastability of a equilibrium points or invariant sets - solving the PDE \eqref{eq:FP} may not be the best approach numerically. For small noise intensities - which are commonly assumed in applications - this is particularly unfortunate since the stochastic dynamics is very close to the zero noise limit $\sigma=0$ on short time scales.

Our approach tries to avoid these difficulties and aims at a natural extension of numerical continuation. We linearize \eqref{eq:simple_ex} around the equilibrium points $x^\pm$ which yields
\be
\label{eq:simple_ex_lin}
dX_t=\left(\mu-3(x^\pm)^2\right)Xdt+\sigma dW_t=:A(x^\pm;\mu)Xdt+\sigma dW_t.
\ee 
Observe that \eqref{eq:simple_ex_lin} is an Ornstein-Uhlenbeck (OU) process \cite{Gardiner}. We will use the variance of the OU process to obtain neighbourhoods of $x^\pm$ within which sample paths of \eqref{eq:simple_ex_lin} stay with high probability. If the initial condition for \eqref{eq:simple_ex_lin} is deterministic then the variance of $X_t$ is 
\be
\label{eq:var_simple_int}
\text{Var}(X_t)=\sigma^2\int_0^t u(t,s)^2 ds
\ee
where $u(t,s)$ is the fundamental solution \cite{Hale} of the system 
\be
\label{eq:lin_hom}
u'=A(x^\pm;\mu)u,\qquad \text{$u(t,t_0)=1$.}
\ee 
Defining $V_t:=\text{Var}(X_t)$ direct differentiation of \eqref{eq:lin_hom} gives that $V_t$ satisfies the ODE
\be
\label{eq:var_simple} 
V'=2A(x^\pm;\mu)V+\sigma^2.
\ee
Since $\mu>0$ we have that $A(x^\pm;\mu)=-2\mu<0$ so that \eqref{eq:var_simple} has a stable equilibrium point at 
\be
\label{eq:var_eq_simple}
\bar{V}(\mu,\sigma)=\frac{-\sigma^2}{2A(x^\pm;\mu)}.
\ee 
Next, consider neighbourhoods of $x^\pm$ given by the variance \eqref{eq:var_eq_simple} of the linearized process (see {e.g.} \cite{BerglundGentz})
\be
\label{eq:Bh}
\cB^\pm(h):=\left\{x\in\R:|x-x^\pm|\leq \sqrt{\bar{V}(\mu,\sigma)}h=\frac{\sigma h}{2\sqrt\mu}\right\}
\ee  
where $\sqrt{\bar{V}}$ can obviously be interpreted as the standard deviation. The main idea of definition \eqref{eq:Bh} is that sample paths of \eqref{eq:simple_ex} stay with high probability inside $\cB^\pm(h)$ if they are started at (or near) $x^\pm$. The parameter $h$ scales the variance neighbourhood and can be used to control the probability to stay inside $\cB^\pm(h)$ for a given time. Hence we can think of $h$ as adjusting the confidence level of our metastable prediction ($h=1$, one standard deviation; $h=2$, two standard deviations; etc.). Figure \ref{fig:fig2}(a) shows $\cB^\pm(h)$ for three different values of $\mu=0.5,1.5,2.5$ with fixed noise $\sigma=0.5$. This demonstrates that \eqref{eq:Bh} can be used to approximate metastability properties for small noise intensities. Obviously all calculations for the SDE \eqref{eq:simple_ex} can be carried out analytically. The open question is whether this approach can be used to construct a general and efficient numerical method. There are several problems that have to be considered:

\begin{enumerate}
 \item[(P1)] \textit{Generalize the construction of $\cB^{\pm}(h)$ to arbitrary $n$-dimensional SDE.} We summarize this well-known construction and the relevant results from probability theory in Section \ref{sec:meta_SDE}.
 \item[(P2)] \textit{Find an efficient way to compute the covariance matrix of an OU-process during numerical continuation and/or for all points on a given equilibrium curve}. The important step to solve this problem efficiently is to realize what information is already available from the deterministic continuation algorithm that can be used to compute the covariance matrix. The main techniques from numerical analysis and control theory are summarized in Section \ref{sec:Lyapunov}. 
 \item[(P3)] \textit{Construct and efficiently compute a test function that detects overlaps of different neighbourhoods $\cB^\pm(h)$.} We suggest a test function based on the distance between ellipsoids. From computational geometry and optimization it is known that the distance can be calculated by solving an optimization problem. The definition of the distance and all computational details are given in Section \ref{sec:ellipsoid}.  
\end{enumerate}

Let us point out again that (P1)-(P3) are essentially all solved (or almost solved) as unconnected problems in various branches of numerical analysis, control theory, dynamical systems, optimization and probability. Our main contribution is to recognize the interplay between the different components which will provide a direct extension of deterministic continuation algorithms to metastable stochastic problems. 

%%%%%%%%%%%%%%%%%%%%%%%%%%%%%%%%%%%%%%%%%%%%%%%%%%%%%
\section{Metastability and Linearization}
\label{sec:meta_SDE}

In this section we address the problem (P1) following Berglund and Gentz \cite{BerglundGentz}. Let $x\in \R^n$ and consider the SDE
\be
\label{eq:SDE}
dx_t=f(x_t;\mu)dt+\sigma F(x_t;\mu) dW_t
\ee
where $W_t=(W_{1,t},W_{2,t},\ldots,W_{k,t})^T$ is standard k-dimensional Brownian motion, $\sigma>0$ controls the noise level, $\mu\in\R$ is a parameter and $f:\R^n\times \R\ra\R^n$ and $F:\R^n\times \R\ra \R^{n\times k}$ are sufficiently smooth maps. Suppose the deterministic part of \eqref{eq:SDE} given by $dx_t=f(x_t;\mu)dt$ has a hyperbolic stable equilibrium point $x^*=x^*(\mu)$ for a given range of parameter values. Using a translation $\bar{x}=x-x^*$ we get
\be
\label{eq:SDE1}
d\bar{x}_t=f(x^*+\bar{x}_t;\mu)dt+\sigma F(x^*+\bar{x}_t;\mu) dW_t.
\ee
Assuming that $F(x^*;\mu)\neq 0$ the approximation of \eqref{eq:SDE1} to lowest order via Taylor expansion is
\be
\label{eq:SDE2}
dX_t=A(x^*;\mu)X_t dt+\sigma F(x^*;\mu)dW_t
\ee
where $A(x;\mu)=(D_xf)(x;\mu)\in\R^{n\times n}$ is the usual Jacobian matrix. Equation \eqref{eq:SDE2} is an $n$-dimensional OU process. We assume that the initial condition $x_0$ is deterministic. The generalization of the variance \eqref{eq:var_simple_int} is the covariance matrix 
\benn
C_t:=\text{Cov}(X_t)=\sigma^2\int_0^t U(t,s)F(x^*;\mu)F(x^*;\mu)^TU(t,s)^Tds
\eenn
where $U(t,s)$ is the fundamental solution of $U'=A(x^*;\mu)U$. Differentiation shows that $C_t$ satisfies the ODE
\be
\label{eq:ODE_cov}
C'=A(x^*;\mu)C+CA(x^*;\mu)^T+\sigma^2F(x^*;\mu)F(x^*;\mu)^T.
\ee
Since $x^*$ is a hyperbolic stable equilibrium point, it follows \cite{Bellman,BerglundGentz} that the eigenvalues of the linear operator
\benn
L(C):=A(x^*;\mu)C+CA(x^*;\mu)^T
\eenn
are given by $\{2\lambda_j\}_{j=1}^n$ where $\lambda_j$ are the eigenvalues of $A(x^*;\mu)$ (and of $A(x^*;\mu)^T$). Therefore \eqref{eq:ODE_cov} has a stable equilibrium solution which is obtained by solving
\be
\label{eq:Lya1}
0=A(x^*;\mu)C+CA(x^*;\mu)^T+\sigma^2F(x^*;\mu)F(x^*;\mu)^T.
\ee
Observe that \eqref{eq:Lya1} is a Lyapunov equation. It is well-known (see e.g. \cite{Jameson}) that the stability of $x^*$ implies the unique solvability of \eqref{eq:Lya1}. For notational simplicity we shall not denote the solution of \eqref{eq:Lya1} as $\bar{C}$ but simply write the symmetric covariance matrix as $C$ or $C(x^*;\mu)$. The main step of solving \eqref{eq:Lya1} numerically at a given parameter value $\mu$ can be found in Section \ref{sec:Lyapunov}. Then one can define a generalization of the variance neighbourhood from Section \ref{sec:ana_ex} as
\be
\label{eq:Bh_big}
\cB(h):=\left\{x\in\R^n:(x-x^*)^TC^{-1}(x-x^*)\leq h^2\right\}
\ee
where $h$ is a parameter that can be interpreted as a probabilistic confidence level. A priori, the set \eqref{eq:Bh_big} may not be well-defined as $C$ may not be invertible. It is well-known from control theory \cite{Sontag,Rugh} that $C$ is invertible if and only if the matrix
\be
\label{eq:controllable}
\text{Con}(A,\sigma F):=[\sigma F \quad \sigma AF \quad \cdots \quad \sigma A^{n-1} F ] \in \R^{n\times nk}\\
\ee
for $A=A(x^*;\mu)$ and $F=F(x^*;\mu)$ has maximal rank; this is sometimes concisely expressed as referring to the matrix pair $(A(x^*;\mu),F(x^*;\mu))$ as controllable \cite{Sontag,Rugh}. From the controllability condition it follows that the invertibility of $C$ is related to the structure of the noise encoded in $F(x^*;\mu)$. In Section \ref{sec:Lyapunov} we discuss the case when $C$ is not invertible. For now assume that $C$ is invertible in which case the set $\cB(h)$ is immediately recognized as a solid ellipsoid with shape matrix 
\benn
Q:=h^2C,\qquad \cB(h)=\left\{x\in\R^n:(x-x^*)^TQ^{-1}(x-x^*)\leq 1\right\}.
\eenn
It can be shown \cite{BerglundGentz} that stochastic sample paths stay in $\cB(h)$ near metastable equilibrium points with high probability. Similar results can also be found in the theory of large deviations \cite{FreidlinWentzell}. It is quite lengthy to state the detailed asymptotic estimates depending on $\sigma$, $h$ and the eigenvalues of $A(x^*;\mu)$. Since we are focusing here on a numerical algorithm we refer the reader to \cite{BerglundGentz} for details. 

%%%%%%%%%%%%%%%%%%%%%%%%%%%%%%%%%%%%%%%%%%%%%%%%%%%%%
\section{The Lyapunov Equation}
\label{sec:Lyapunov}

The next step is the numerical solution of the Lyapunov equation for a given metastable equilibrium $x^*(\mu)$ as well as for an entire branch of equilibrium points obtained via continuation $\gamma=\{(x^*(\mu),\mu)\}\subset \R^n\times \R$. The algebraic equation \eqref{eq:Lya1} is a uniquely solvable Lyapunov equation of the form
\be
\label{eq:Lya_num}
AC+CA^T+B=0
\ee
where we are going to use the shorthand notations $A=A(x^*;\mu)$ and $B:=\sigma^2F(x^*;\mu)F(x^*;\mu)^T$ from now on. Lyapunov equations have been studied in various branches of mathematics \cite{GajicQureshi}. Recall \cite{Jameson} that if one sorts the elements of $C$ and $B$ in vector form 
\benn
c^T=(c_{11},c_{21},\cdots,c_{12},\cdots)^T \qquad \text{and}\qquad b^T=(b_{11},b_{21},\cdots,b_{12},\cdots)^T
\eenn
then \eqref{eq:Lya_num} can be rewritten as a standard linear system
\be
\label{eq:big_lin_Lya}
[I \otimes A +A \otimes I]c=-b
\ee
where $\otimes$ denotes the Kronecker product \cite{GolubVanLoan}. The problem of efficient numerical solution of \eqref{eq:Lya_num} or \eqref{eq:big_lin_Lya} (and of several generalizations) has attracted considerable attention in numerical analysis and control theory \cite{GajicQureshi}. For our situation several new aspects arise since we want to solve \eqref{eq:Lya_num} along an entire equilibrium branch $\gamma$:

\begin{enumerate}
 \item All standard numerical continuation algorithms require an approximation of the $n\times (n+1)$ Jacobian matrix $(D_{(x,\mu)}f)(x^*(\mu_1),\mu_1)$ to compute a point $(x^*(\mu_2),\mu_2)\in\gamma$ starting from $(x^*(\mu_1),\mu_1)\in \gamma$. Therefore, the matrix $A=(D_xf)(x^*(\mu_1);\mu_1)$ is available at each continuation step. Furthermore, computing the matrix $B$ requires at most one matrix multiplication at a given point $(x^*(\mu_1),\mu_1)$; for purely additive noise $B$ can even be precomputed for all equilibrium points.
 \item Solving \eqref{eq:Lya_num} at $(x^*(\mu_1),\mu_1)$ gives a matrix $C(x^*(\mu_1);\mu_1)$. If $|\mu_1-\mu_2|$ is small then $C(x^*(\mu_1);\mu_1)$ is already an excellent initial guess to find $C(x^*(\mu_2);\mu_2)$! Hence, except for the first point on the equilibrium curve, we always have an initial guess available for iterative methods.  
\end{enumerate}

The observations suggest that computing the covariance $C$ should be relatively easy. We decided to focus on three different approaches which we briefly review here. Due to a good initial guess, the most natural choice are iterative methods. Consider the reformulation \eqref{eq:big_lin_Lya} and define $\cA:=[I \otimes A +A \otimes I]$. Then the standard Gauss-Seidel iteration \cite{StoerBulirsch} is given by
\be
\label{eq:Gauss_Seidel}
c^{(k+1)}=c^{(k)}-M_{GS}^{-1}(\cA c^{(k)}+b)
\ee
where $M_{GS}$ is the matrix obtained from $\cA$ by setting all entries above the diagonal ($\cA_{ij}$, $j>i$) to zero. The iteration is terminated when $\|c^{(k+1)}-c^{(k)}\|<tol$ where $tol$ is a given tolerance. Other possibilities for iterative methods include the Jacobi method and successive overrelaxation (SOR) methods \cite{StoerBulirsch}. For large sparse Lyapunov equations several special methods have been suggested including alternating-direction-implicit (ADI) by Wachspress \cite{Wachspress} and special SOR methods by Starke \cite{Starke}. We shall not consider the special methods here although they should definitely be relevant for large scale bifurcation problems \cite{Salingeretal}.

Another well-known method for the iterative solution of \eqref{eq:Lya_num} is Smith's algorithm \cite{Smith}. The first step is to fix a scalar $q>0$ and consider the matrices 
\beann
K&:=&2q(qI-A)^{-1}B(qI-A^T)^{-1},\\
G&:=&(qI-A)^{-1}(qI+A).
\eeann
Direct matrix multiplication shows that \eqref{eq:Lya_num} is equivalent to solving
\be
\label{eq:Lya_iter}
C=K+GCG^T.
\ee
The iteration of \eqref{eq:Lya_iter} converges linearly. Smith observed that with initial guess $C^{(0)}=K$ the iteration
\be
\label{eq:Smith}
C^{(k+1)}=C^{(k)}+G^{2^k}C^{(k)}\left(G^{2^k}\right)^T,\qquad k=0,1,2,\ldots
\ee 
obtained by squaring $G$ at each step converges quadratically. The algorithm is terminated when $\|C^{(k+1)}-C^{(k)}\|<tol$. Using ADI theory the optimal $q>0$ can be found and the error has been calculated \cite{Wachspress,GajicQureshi}; we will simply fix $q=0.1$ which is the classical choice by Smith \cite{Smith1}. Observe that Smith's algorithm does not use an initial guess. 

There are also several direct (non-iterative) algorithms available. The most important techniques were suggested in the 1970s \cite{BartelsStewart,BelangerMcGillivray,GolubNashvanLoan} and have become standard methods for the numerical solution of \eqref{eq:Lya_num}. The Bartels-Stewart algorithm \cite{BartelsStewart} requires to compute the real Schur decomposition of $A$ given by $U^TAU=R$ where $U\in\R^{n\times n}$ is orthogonal and $R\in\R^{n\times n}$ is upper quasi-triangular ({i.e.} diagonal with possible $2\times 2$ blocks on the diagonal corresponding to complex eigenvalues). Then \eqref{eq:Lya_num} can be transformed to
\be
\label{eq:quasi_up}
RY+YR^T=W
\ee
where $W=-U^TBU$ and $Y=U^TCU$. The right-hand side $W$ can be obtained by solving $UW=BU$ for $W$. Solving \eqref{eq:quasi_up} requires the solution of an upper quasi-triangular system which is straightforward. Then one can solve $CU=UY$ for $C$ to get the final result. The Bartels-Stewart algorithm can also be helpful for our problem as it can be used to solve the problem at the first continuation point and it applies when Gauss-Seidel iteration fails as a ``fall-back'' strategy.

In Section \ref{sec:num_ex1} we are going to compare the performance of the Bartels-Stewart algorithm, Smith's method and Gauss-Seidel iteration for a practical numerical continuation problem. Once we have the covariance matrix $C$ it is important to check whether $C^{-1}$ exists so that \eqref{eq:Bh} is a well-defined ellipsoid. We are going to illustrate why such a test is important. In control theory \cite{KurzhanskiiValyi} it is well-known how to define ellipsoids in a degenerate case when the shape matrix $Q=h^2C$ is only positive semidefinite (see also Section \ref{sec:ellipsoid})
\be
\label{eq:Bh_gen}
\cB(h):=\left\{x\in\R^n:v^Tx\leq v^Tx^*+(v^TQv)^{1/2}\quad \forall v\in\R^n\right\}.
\ee
Now consider the example
\benn
A=A(x^*;\mu):=\left(
\begin{array}{cc}
-2 & 0 \\ 0 & -1 \\
\end{array}
\right),\qquad 
B=\sigma^2F(x^*;\mu)F(x^*;\mu)^T:=\left(
\begin{array}{cc}
\sigma^2 & 0 \\ 0 & 0 \\
\end{array}
\right)
\eenn
which corresponds to a stable hyperbolic equilibrium in $\R^2$ with additive noise on the first component only. Even without solving for $C$ we can compute the matrix \eqref{eq:controllable}
\benn
\text{Con}(A(x^*;\mu),\sigma F(x^*;\mu))=\left(
\begin{array}{cccc}
\sigma & 0 & -2\sigma & 0 \\ 0 & 0 & 0 & 0\\
\end{array}
\right)
\eenn 
so that $C$ is not invertible since $\text{rank}(\text{Con}(A(x^*;\mu),\sigma F(x^*;\mu)))=1$. Indeed, we easily find that solving the Lyapunov equation gives
\benn
C=\left(
\begin{array}{cc}
\sigma^2/4 & 0 \\ 0 & 0 \\
\end{array}
\right)\qquad \Rightarrow \quad
Q=\left(
\begin{array}{cc}
(h\sigma)^2/4 & 0 \\ 0 & 0 \\
\end{array}
\right)
\eenn
Assuming for simplicity that $x^*=(0,0)^T$ we get that the set $\cB(h)$ defined in \eqref{eq:Bh_gen} is given by 
\beann
\cB(h)&=&\left\{(x_1,x_2)\in\R^2:v_1x_1+v_2x_2\leq \frac{\sqrt{v_1^2\sigma^2 h^2}}{2}\quad \forall (v_1,v_2)\in\R^2\right\}\\
&=&\left\{(x_1,x_2)\in\R^2:|x_1|\leq \sigma h/2,~x_2=0\right\}.
\eeann
The ellipsoid $\cB(h)$ is a degenerate interval which reflects that the degenerate noise terms only act on the $x_1$-coordinate. Detecting such a degenerate (or near-degenerate) noise is clearly important in applications as this identifies directions along which metastable escapes are unlikely. A simple test for this degeneracy is to compute the singular value decomposition (SVD) \cite{GolubVanLoan} of $C$.

%%%%%%%%%%%%%%%%%%%%%%%%%%%%%%%%%%%%%%%%%%%%%%%%%%%%%
\section{Ellipsoids and the Testfunction}
\label{sec:ellipsoid}

Suppose we have two covariance matrices $C_{1,2}$ for given set of parameter values at $(\mu,x^*_{1,2}(\mu))=(\mu,x^*_{1,2})$ and we are interested in detecting the distance between the associated ellipsoids as large/small distances are expected to correspond to long/short travel times of sample paths. Denote the shape matrices of the ellipsoids by $Q_{1,2}=h^2C_{1,2}$ and a general ellipsoid by 
\benn
\cE=\cE(x^*,Q)=\{x\in\R^n:(x-x^*)^TQ^{-1}(x-x^*)\leq 1\}.
\eenn
The idea of considering covariance ellipsoids and their overlaps is not new. For example, the idea is used in satellite tracking for collision avoidance  \cite{Alarcon-RodriguezMartinez-FadriqueKlinkrad}. In computational geometry and robotics one often considers the minimum-volume enclosing ellipsoid of an object, also called the L\"{o}wner-John ellipsoid \cite{RimonBoyd}. In statistics an ellipsoidal distance defined via the covariance matrix, the so-called Mahalanobis distance \cite{Mahalanobis}, is often used \cite{Rosenbaum}. Various method have been proposed to detect ellipsoid overlaps ranging from Gr\"{o}bner bases \cite{Buchberger}, analytical representation formulas \cite{AlfanoGreer}, reformulation as an eigenvalue problem \cite{RimonBoyd}, local approximation by balls \cite{LinHan} to polyhedral approximations \cite{GilbertJohnsonKeerthi,Bobrow}. Here we will adapt an idea based on calculating the distance between ellipsoid by solving an optimization problem which has several advantages to be discussed below. The support function of an ellipsoid is \cite{KurzhanskiiValyi}
\benn
\sigma_{\cE}(v):=\sup_{x\in \cE}v^Tx=v^Tx^*+(v^TQv)^{1/2}.
\eenn
The Hahn-Banach Theorem \cite{AuslenderTeboulle} gives that an ellipsoid with a positive semi-definite shape matrix can be defined as 
\be
\label{eq:ell_gen}
\cE:=\left\{x\in\R^n:v^Tx\leq v^Tx^*+(v^TQv)^{1/2}\quad \forall v\in\R^n\right\}.
\ee
A measure of the distance between two ellipsoids \cite{Kurzhanskiy} is given by
\bea
\label{eq:ell_dist}
\delta=\delta(\cE(x^*_1,Q_1),\cE(x^*_2,Q_2))&=&\max_{\|v\|=1}\left(-\sigma_{\cE_1}(-v)-\sigma_{\cE_2}(v)\right)\nonumber\\
&=&\max_{\|v\|=1}\left(v^Tx_1^*-(v^TQ_1v)^{1/2}-v^Tx_2^*-(v^TQ_2v)^{1/2}\right).
\eea
The distance $\delta$ and its definition have several advantages for detecting metastability. The definition also applies immediately if the matrices $Q_{1,2}$ are degenerate. For example, if we are interested in the distance of a covariance ellipsoid $\cE(x^*_1,Q_1)$ to an unstable equilibrium point $x^*_2$ ({e.g.} the saddle point in Section \ref{sec:ana_ex}) we can just set $Q_2=0$ and still consider the distance $\delta$. We can even replace the ellipsoid with a more general convex set $\cH$ if the support function $\sigma_{\cH}(v)$ is easy to calculate. The main advantage is that $\delta$ is also a test function since

\begin{itemize}
 \item $\delta(\cE(x^*_1,Q_1),\cE(x^*_2,Q_2))>0$ if the two ellipsoids are disjoint,
 \item $\delta(\cE(x^*_1,Q_1),\cE(x^*_2,Q_2))=0$ if the ellipsoids touch at a point, and
 \item $\delta(\cE(x^*_1,Q_1),\cE(x^*_2,Q_2))<0$ if the ellipsoids intersect.
\end{itemize}

%A simple illustration of the function \eqref{eq:ell_dist} is shown in Figure \ref{fig:fig3} for the two ellipsoids 
%\be
%\label{eq:test_dist}
%(x_1^*,Q_1)=\left(
%\left(\begin{array}{c} 0 \\ 0\\\end{array}\right),
%\left(\begin{array}{cc} 3 & 0 \\ 0 & 1\\\end{array}\right)
%\right)\qquad \text{and} \qquad 
%(x_2^*,Q_2)=\left(
%\left(\begin{array}{c} \rho \\ 0\\\end{array}\right),
%\left(\begin{array}{cc} 1 & 0 \\ 0 & 2\\\end{array}\right)
%\right)
%\ee
%for $\rho=4$ and $\rho=2$. For the first case, the two ellipsoids are disjoint with distance $\delta_{\rho=4}\approx 1.2679$ while the overlap in the second case with $\delta_{\rho=2}\approx -0.7321$.\\
%
%\begin{figure}[htbp]
%\psfrag{a}{(a)}
%\psfrag{b}{(b)}
%\psfrag{E1}{$\cE_1$}
%\psfrag{E2}{$\cE_2$}
%\psfrag{x1}{$x_1$}
%\psfrag{x2}{$x_2$}
%	\centering
%		\includegraphics[width=1\textwidth]{./fig03.eps}
%	\caption{\label{fig:fig3}Simple illustration for the test function \eqref{eq:ell_dist} for the two ellipsoids given in \eqref{eq:test_dist}. The points on $\partial \cE_{1,2}$ that realize the distance \eqref{eq:ell_dist} are connected by a (thick) line. (a) $\rho=4$: The distance $\delta\approx 1.2679$ is positive so that no overlap occurs. (b) $\rho=2$: The distance $\delta\approx -0.7321$ is negative signalling overlap.}
%\end{figure} 

Therefore the distance \eqref{eq:ell_dist} is a test (or bifurcation) function if we want to check how likely metastable transitions occur in our SDE \eqref{eq:SDE}. Note that the precise number of noise-induced transitions cannot be inferred from $\delta$ as we have not made any assumptions about global dynamics; but see Section \ref{sec:Kramer}. Observe that \eqref{eq:ell_dist} is a classical nonlinear optimization (or nonlinear programming) problem. In standard minimization form with a differentiable constraint it can be written as
\be
\label{eq:ell_dist_opt}
\left\{\begin{array}{rl}
\min&\left(-v^Tx_1^*+(v^TQ_1v)^{1/2}+v^Tx_2^*+(v^TQ_2v)^{1/2}\right)=:\min ~G(v) ,\\
\text{subject to}& 0=\|v\|^2-1=:g(v).\\
\end{array}\right.
\ee
and we obtain a solution to \eqref{eq:ell_dist} by the negative solution value of \eqref{eq:ell_dist_opt}. Many efficient algorithms for the solution of \eqref{eq:ell_dist_opt} are available \cite{NocedalWright}. In particular, many iterative schemes are known among which sequential quadratic programming (SQP) \cite{Han,Powell} has turned out to be among the most powerful techniques. Here we simply use this approach which solves a quadratic programming problem at iteration step $k$ given by 
\be
\label{eq:ell_dist_opt1}
\left\{\begin{array}{rl}
\min &\left(\frac12 w^TH_kw+\nabla G(v_k)^Tw\right),\\
\text{subject to}& \nabla g(v_k)^Tw+g(v_k)=0,\\
\end{array}\right.
\ee
where $H_k=(\nabla^2_{v}L)(v_k,u_k)$ is the Hessian of the Lagrangian $L(v,u):=G(v)-u^Tg(v)$ and $u_k\in\R$ is an approximation of the Lagrange multiplier. If the solution of \eqref{eq:ell_dist_opt1} at step $k$ is denoted by $w_k$ then the main iteration step is $v_{k+1}=v_k+\alpha_kw_k$ for a given step length $\alpha_k>0$. It is important to note that $\nabla G(v_k)$, $\nabla g(v_k)$ and $H_k$ can be supplied in explicit form
\beann
\nabla G(v_k)&=&-x_1^*+\frac12(v_k^TQ_1v_k)^{-1/2}(Q_1+Q_1^T)x_k+x_2^*+\frac12(v_k^TQ_2v_k)^{-1/2}(Q_2+Q_2^T)x_k,\\
\nabla g(v_k)&=&2x_k,\\
\{(H_k)_{ij}\}_{i,j=1}^n&=&-\frac14(v_k^TQ_1v_k)^{-3/2}[(Q_1+Q_1^T)v_k]_i[(Q_1+Q_1^T)v_k]_j+\frac12(v_k^TQ_1v_k)[Q_1+Q_1^T]_{ij}\\
&&-\frac14(v_k^TQ_2v_k)^{-3/2}[(Q_2+Q_2^T)v_k]_i[(Q_2+Q_2^T)v_k]_j+\frac12(v_k^TQ_2v_k)[Q_2+Q_2^T]_{ij}
\eeann
which avoids the computation of finite difference approximations during the optimization iteration. We use a standard quasi-Newton line-search method to solve \eqref{eq:ell_dist_opt1}. The iterative algorithm stops when the solution, solution values and constraints are below a given tolerance. The details of this part of the algorithm will not be discussed here and details can be found in \cite{NocedalWright,MatLab2010b}.

As for the Gauss-Seidel method, it is very important to point out that the iterative solution of \eqref{eq:ell_dist} can be used efficiently during continuation. Given a fixed point $x^*(\mu_1)$ at parameter values $\mu_1$ we obtain a solution $v(\mu_1)$ to \eqref{eq:ell_dist} by solving \eqref{eq:ell_dist_opt}. For an equilibrium point continuation step from $\mu_1$ to $\mu_2$ we have that $|\mu_1-\mu_2|$ is small so that $v(\mu_1)$ can be used as a very good initial guess for the optimization problem to be solved with parameter values $\mu_2$.

%%%%%%%%%%%%%%%%%%%%%%%%%%%%%%%%%%%%%%%%%%%%%%%%%%%%%
\section{Algorithm Summary}
\label{sec:main}

In this section, we summarize the main steps of our algorithm which augments deterministic numerical continuation. Consider the SDE
\be
\label{eq:SDE_final}
dx_t=f(x_t;\mu)dt+\sigma F(x_t;\mu)dW_t,\qquad \text{for $x_t\in\R^n$ and $\mu\in\R$.}
\ee
We assume that a stable equilibrium $x^*(\mu_0)$ for the deterministic part of \eqref{eq:SDE_final} is given (or it can be found e.g.~using Newton's method \cite{StoerBulirsch}) so that $f(x^*(\mu_0);\mu)=0$. Then define 
\benn
A(x^*(\mu_0);\mu_0)=(D_xf)(x^*(\mu_0);\mu_0).
\eenn
Using the Bartels-Stewart algorithm (see Section \ref{sec:Lyapunov}) we solve
\benn
0=A(x^*(\mu_0);\mu_0)C+CA(x^*(\mu_0);\mu_0)^T+\sigma^2F(x^*(\mu_0);\mu_0)F(x^*(\mu_0);\mu_0)^T
\eenn
for the covariance matrix $C=C(x^*(\mu_0);\mu_0)$. This completes the initialization step. The main iterative step of the algorithm is as follows:

\begin{itemize}
 \item[(A1)] Choose a step length $\beta_k$ and set $\mu_{k}=\mu_{k-1}+\beta_k$. Solve the continuation problem $f(x;\mu_{k})=0$ for the new equilibrium $x^*(\mu_{k})$ with starting point $x^*(\mu_{k-1})$ (see e.g. \cite{Kuznetsov}).
 \item[(A2)] Consider the Lyapunov equation 
 \benn
 0=A(x^*(\mu_{k});\mu_{k})C+CA(x^*(\mu_{k});\mu_{k})^T+\sigma^2F(x^*(\mu_{k});\mu_{k})F(x^*(\mu_{k});\mu_{k})^T
 \eenn
  and solve it for $C$, preferably using an iterative algorithm with initial guess $C(x^*(\mu_{k-1});\mu_{k-1})$. This yields the new covariance matrix $C(x^*(\mu_{k});\mu_{k})$. Define the shape matrix $Q(x^*(\mu_{k});\mu_{k}):=h^2C(x^*(\mu_{k});\mu_{k})$ for a given confidence level $h$.
 \item[(A3)] Given two shape matrices $Q_{1,k}:=Q(x^*_1(\mu_{k});\mu_{k})$ and $Q_{2,k}:=Q(x^*_2(\mu_{k});\mu_{k})$ for different stable equilibria, compute the distance $\delta=\delta(\mu_k,h)$ between the two ellipsoids $\cE_1(x_1^*(\mu_{k});Q_{1,k})$ and $\cE_2(x_2^*(\mu_{k});Q_{2,k})$ solving the optimization problem \eqref{eq:ell_dist} using an iterative method such as SQP (see Section \ref{sec:ellipsoid}) with initial conditions obtained from the iteration step $k-1$.
\end{itemize}

As an output we get the following parameterized families
\begin{itemize}
 \item equilibrium points $\{(x^*(\mu);\mu)\}_\mu$ from numerical continuation,
 \item ellipsoids $\{\cE(x^*(\mu),C(x^*;\mu))\}_\mu$ from solving Lyapunov equations, and
 \item mutual distances $\{\delta(\mu,h)\}_\mu$ from solving a nonlinear programming problems. 
\end{itemize}

The ellipsoids $\{\cE(x^*(\mu),C(x^*;\mu))\}_\mu$ provide locally rigorous estimates for metastability \cite{BerglundGentz}. The distance $\delta(\mu_k,h)$ between two ellipsoids $Q_{1,k}$ and $Q_{2,k}$ gives an indicator for global transitions occurring from $Q_{1,k}$ to $Q_{2,k}$ or vice versa, based on the assumption that larger distances correspond to lower switching probabilities. Section \ref{sec:num_ex1} shows that using $\delta(\mu_k,h)$ works nicely in practice. Nevertheless, it may be desirable to obtain rigorous estimates if global assumptions are made; Section \ref{sec:Kramer} augments the algorithm in this direction for a special case. 

Furthermore, it looks intuitive to consider higher-order moments of the fully nonlinear stochastic process described by the SDE \eqref{eq:SDE}. However, the ODEs for higher-order moments usually do not form a closed system \cite{Gardiner} such as \eqref{eq:SDE2}. Observe carefully that if a set of moment equations forms a finite-dimensional closed system (or can approximately be closed \cite{Socha}) then a modified version of steps (A1)-(A2) should carry over to this situation since equilibria for the moment ODEs satisfy an algebraic equation which can again be solved iteratively with a good initial guess from the previous continuation step.

We note that the algorithmic steps (A2)-(A3) can be used as a post-processing tool for an existing numerical continuation curve $(x^*(\mu);\mu)$. This is not as efficient as combining (A1)-(A3) as it requires re-building the matrices $A(x^*;\mu)$.
In summary, we have obtained local approximate information about a system of stochastic differential equations using a completely deterministic continuation algorithm. The additional computations required to obtain this information are easy to implement in a classical continuation algorithm and/or bifurcation software package. The iterative solution procedures for the Lyapunov equation and the ellipsoid distances are expected to make the algorithm computationally very efficient.

%%%%%%%%%%%%%%%%%%%%%%%%%%%%%%%%%%%%%%%%%%%%%%%%%%%%%
\section{Neural Competition and Bistability}
\label{sec:num_ex1}

In this section we are going to test our algorithm for the situation where the deterministic dynamical system has two stable coexisting equilibrium points (``bistability''). The differential equations we are going to study are based on ODEs modelling a two-cell inhibitory neural network \cite{CurtuShpiroRubinRinzel,Curtu}. The goal is to describe competition between two neural populations. For example, such a situation can occur due to ambiguous external stimuli \cite{LaingChow} inducing a bistable behaviour in the neuronal system. A typical example is binocular rivalry \cite{EinhaeuserMartinKoenig} where switching between different visual perceptions occurs. This situation can be modelled \cite{ShpiroCurtuRinzelRubin} by considering the (fast) spatially averaged firing rates $x_{1,2}$ of two neural populations and two associated (slow) time fatigue accumulation variables $y_{1,2}$. The resulting ODEs are 
\be
\label{eq:Curtu_full}
\begin{array}{rcl}
x_1' & = & -x_1+S(I_c-\beta x_2-gy_1),\\
x_2' & = & -x_2+S(I_c-\beta x_1-gy_2),\\
y_1' & = & \epsilon(x_1-y_1),\\
y_2' & = & \epsilon(x_2-y_2),\\
\end{array}
\ee
where $I_c$ is the main bifurcation parameter and the sigmoid-shaped gain function $S:\R\ra\R$ is often chosen \cite{CurtuShpiroRubinRinzel} in numerical simulations and continuation calculations as
\benn
S(u):=\frac{1}{1+\exp(-r(u-\theta))}.
\eenn
We adopt this choice and also fix the parameters
\be
\label{eq:Curtu_fixed_paras}
\beta=1.1,\qquad g=0.5,\qquad r=10,\qquad\theta=0.2
\ee
so that our calculations are a direct extension of numerical continuation in \cite{Curtu}. The parameter $0\leq \epsilon\ll 1$ describes the time scale separation between the fast and slow variables. We are only going to consider \eqref{eq:Curtu_full} in the singular limit $\epsilon=0$ of perfect time scale separation. The equations 
\be
\label{eq:Curtu_fast}
\begin{array}{rcl}
x_1' & = & -x_1+S(I_c-\beta x_2-gy_1)=:f_1(x),\\
x_2' & = & -x_2+S(I_c-\beta x_1-gy_2)=:f_2(x),\\
\end{array}
\ee
are also called the fast subsystem of \eqref{eq:Curtu_full} where $y_{1,2}$ are regarded as parameters. For an introduction to the theory of fast-slow systems and singular limits see \cite{Jones,MisRoz}; an example how fast subsystem bifurcation analysis can form a building block of bifurcation analysis for the case $\epsilon>0$ can be found in \cite{GuckenheimerKuehn1,GuckenheimerKuehn3}. Since \eqref{eq:Curtu_fast} is a model for the activity of (finite) neuronal populations there are various natural stochastic effects such as channel noise \cite{Fox}, input noise \cite{Tuckwell}, neuronal background noise \cite{Koch} and external noise in experiments/observations \cite{KassVenturaBrown}. Therefore it is reasonable to extend \eqref{eq:Curtu_fast} to the SDE
\be
\label{eq:Curtu_SDE}
\left(\begin{array}{c}
dx_1 \\
dx_2 \\
\end{array}\right)
=
\left(\begin{array}{c}
-x_1+S(I_c-\beta x_2-gy_1)\\
-x_2+S(I_c-\beta x_1-gy_2)\\
\end{array}\right)dt
+\sigma^2F(x)dW_t
\ee
where $F:\R^2\ra \R^{2\times 2}$. Furthermore, we fix the slow variables to
\be
\label{eq:Curtu_fixed_paras1}
y_1=0.7 \qquad \text{and} \qquad y_2=0.75
\ee
which introduces a slight asymmetry into the system. Both slow variables also lie within plausible ranges as considered in \cite{CurtuShpiroRubinRinzel}.

\begin{figure}[htbp]
\psfrag{a}{(a)}
\psfrag{b}{(b)}
\psfrag{c}{(c)}
\psfrag{d}{(d)}
\psfrag{e}{(e)}
\psfrag{f}{(f)}
\psfrag{MFPT}{$T_p$}
\psfrag{Ic}{$I_c$}
\psfrag{t}{$t$}
\psfrag{del}{$\delta$}
\psfrag{x1}{$x_1$}
\psfrag{x2}{$x_2$}
	\centering
		\includegraphics[width=1\textwidth]{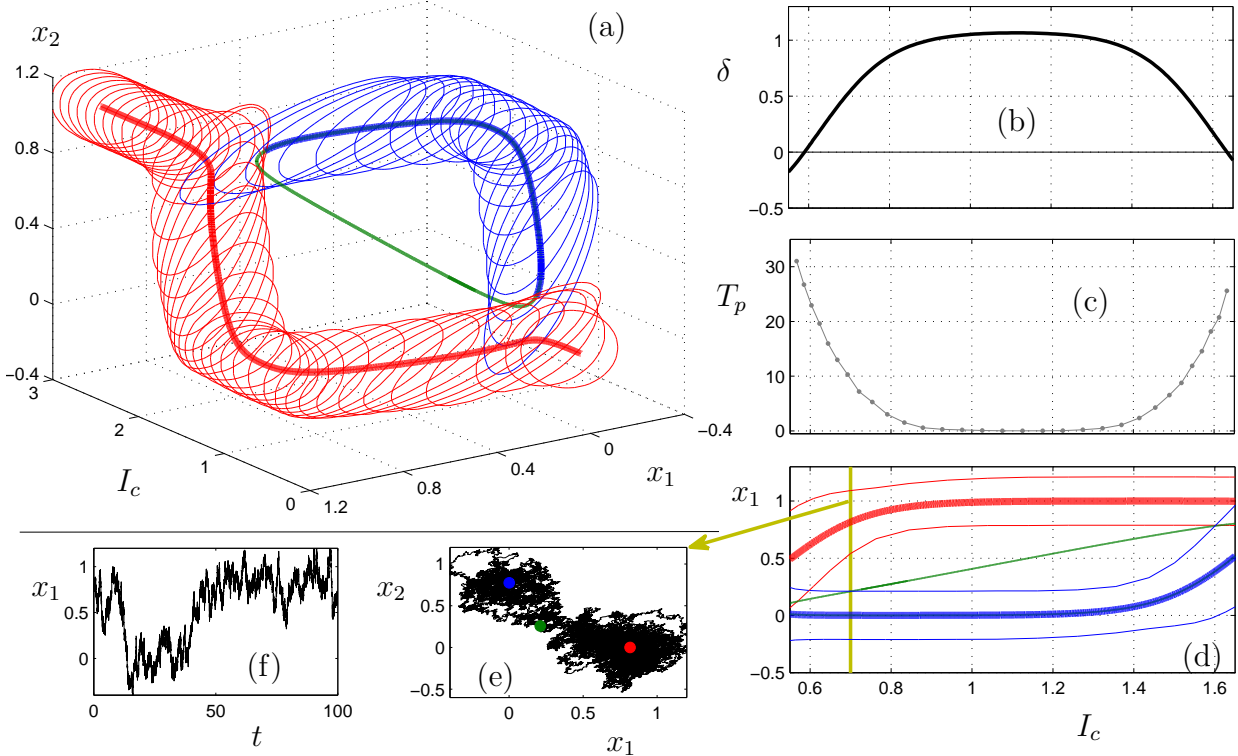}
	\caption{\label{fig:fig4}Continuation results for \eqref{eq:Curtu_SDE} with parameter values \eqref{eq:Curtu_fixed_paras} and \eqref{eq:Curtu_fixed_paras1}. The noise terms are given in \eqref{eq:Curtu_noise}. (a) The thick curves (red and blue) show stable equilibrium point branches continued in the main bifurcation parameter $I_c$. There are two saddle node bifurcations on the isolated bifurcation curve (isola) from a stable node (blue, thick curve) to a saddle (green, thin curve). We also show some of the two-dimensional ellipsoids $\cE_{1,2}$ calculated in $(x_1,x_2)$ phase space for fixed parameter values $I_c$ and embedded in $(x_1,x_2,I_c)$-coordinates; here we use $h=1$ as the confidence parameter in the definition of the covariance matrix. (b) Calculation of the distance $\delta=\delta(\cE_1,\cE_2)$ between the ellipsoids in the bistable regime as defined in \eqref{eq:ell_dist}. (c) Mean number of passages $T_p$ of a trajectory between the two stable equilibria over a time interval $[0,1000]$; averaged results over 100 sample paths obtained via direct numerical integration of the SDE are shown (grey, dots indicate grid in $I_c$). (d) Projection of (a) onto $(I_c,x_1)$ where the $x_1$ maxima and minima of the ellipsoids have been connected to form tubes around the stable equilibrium branches. (e)-(f) Direct numerical SDE simulation for $I_c=0.7$ (as indicated by the arrow from (d)). The colored dots correspond to the equilibria.}
\end{figure} 

Figure \ref{fig:fig4} shows a continuation calculation for the neuronal competition model \eqref{eq:Curtu_SDE} with parameter values \eqref{eq:Curtu_fixed_paras} and \eqref{eq:Curtu_fixed_paras1}. The additive noise terms are given by
\be
\label{eq:Curtu_noise}
\sigma^2F(x^*)F(x^*)^T=\sigma^2
\left(\begin{array}{cc}
1 & 0.4 \\ 0.4 & 1\\
\end{array}\right)\qquad \text{for $\sigma=0.3$.}
\ee
The deterministic equilibrium continuation has been carried out using the Moore-Penrose algorithm \cite{Kuznetsov,MatCont} with fixed continuation step size $0.001$. For the the computation of the covariance ellipsoids and the distance between them we refer to the summary of our algorithm in Section \ref{sec:main}. Figures \ref{fig:fig4}(a) and (d) visualize the ellipsoids and Figure \ref{fig:fig4}(b) shows the distance between the ellipsoids defined by \eqref{eq:ell_dist}. We find two regions where the distance is negative and overlaps between ellipsoids occur. Hence we expect that equilibrium points in the parameter regions with overlaps are only weakly metastable and relatively frequent noise-induced switching between different neuronal activity patterns occurs. This conjecture is confirmed in Figure \ref{fig:fig4}(c) where the mean number of noise-induced passages $T_p$ between two stable equilibrium points $p_{1,2}$ is shown during a fixed time interval; more precisely, consider a trajectory $\gamma(t)$, fix some small $\rho>0$ and define
\benn
\begin{array}{lcl}
T^{1\ra 2}_p(\gamma)&:=&\#\left\{(t_1,t_2):t_1<t_2<T_{max},\|\gamma(t_1)-p_1\|<\rho,t_2=\inf\{t:t>t_1,\|\gamma(t)-p_2\|<\rho\}\right\},\\
T^{2\ra 1}_p(\gamma)&:=&\#\left\{(t_1,t_2):t_1<t_2<T_{max},\|\gamma(t_1)-p_2\|<\rho,t_2=\inf\{t:t>t_1,\|\gamma(t)-p_1\|<\rho\}\right\},\\
\end{array}
\eenn
which just count the number of times a trajectory starting from a small ball near $p_1$ reaches as small ball near $p_2$ and vice versa. Then we can average the results over different realizations of the noise (i.e. over different paths $\gamma$)
\be
\label{eq:Tp}
T_p:=\E[T^{1\ra 2}_p+T^{2\ra 1}_p].
\ee
For Figure \ref{fig:fig4}(c) the parameters $\rho=0.05$ and $T_{max}=1000$ have been used and the expected value in \eqref{eq:Tp} has been computed over 100 sample paths. Note carefully that distance function in Figure \ref{fig:fig4}(b) predicts the qualitative shape of the passage time distribution $T_p$ very nicely.

\begin{figure}[htbp]
\psfrag{a}{\scriptsize{(a)}}
\psfrag{b}{\scriptsize{(b)}}
\psfrag{steps}{\scriptsize{iter}}
\psfrag{time}{\scriptsize{time}}
\psfrag{tol}{\scriptsize{tol}}
	\centering
		\includegraphics[width=1\textwidth]{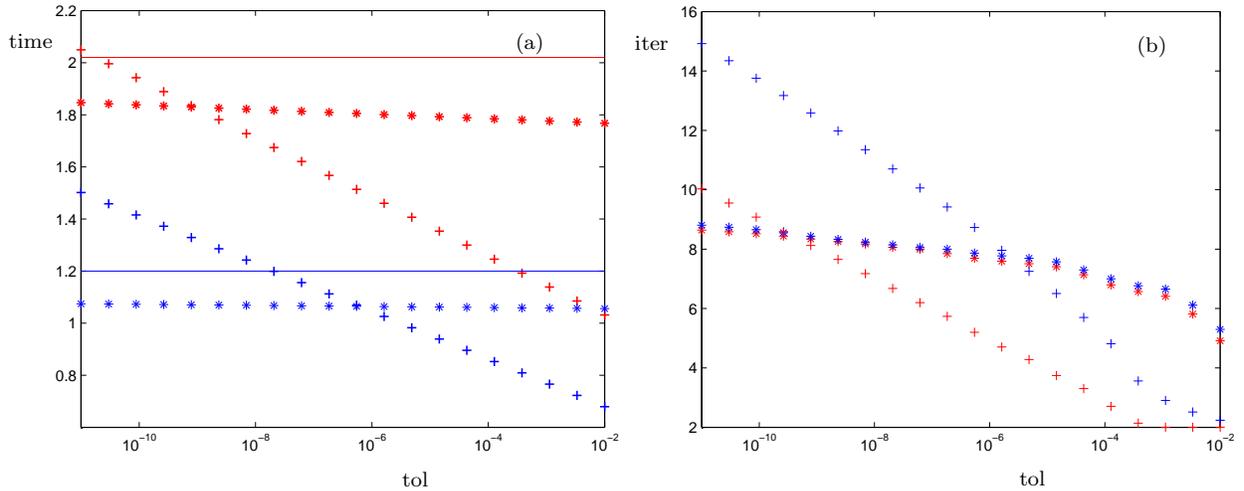}
	\caption{\label{fig:fig5}Computation of the covariance matrices along the entire two stable equilibrium continuation curves (see Figure \ref{fig:fig4}) for given tolerances (logarithmic abscissa). The colors (red/blue) indicate the equilibrium branch as in Figure \ref{fig:fig4}. (a) Total computation time (in seconds, linear interpolation of times is shown) for different tolerances of $\|C^{(k)}-C^{(k-1)}\|\leq tol$. Gauss-Seidel iteration (crosses), Smith iteration (stars) and a direct solution via the Bartels-Stuart algorithm (as implemented in \cite{MatLab2010b}) are compared. (b) Comparison of the total number of iteration steps for the Gauss-Seidel (crosses) and Smith (stars) algorithms are compared.}
\end{figure} 

We shall not investigate the dynamical implications from our method here but focus on the performance of the algorithm. As a starting point we use the two continuation curves of stable equilibrium points shown in Figure \ref{fig:fig4}. For each curve we calculate the covariance matrix by solving the Lyapunov equation for each point on the continuation curve using Gauss-Seidel and Smith iterations as well as the Bartels-Stewart algorithm. For the Gauss-Seidel algorithm we use as the starting point of the iteration the covariance matrix from the previous point on the equilibrium curve. Figure \ref{fig:fig5} shows the computation time as well as the average number of iteration steps along the equilibrium curve for different tolerances of the iteration termination condition
\benn
\|C^{(k)}-C^{(k-1)}\|\leq tol.
\eenn 
We see that for relatively low tolerances between $10^{-2}$ and $10^{-7}$ the iterative solution using the Gauss-Seidel method seems to perform best. This is not surprising since it is the only method that uses the previous point on the curve of equilibria which is expected to be an excellent initial guess. For higher tolerances and high-precision computation Smith's algorithm as well as the exact Bartels-Stewart method seem to be preferable. Since Smith's algorithm always converges quadratically this is again expected in comparison to Gauss-Seidel. Using SOR or ADI iterative techniques or considering larger systems could potentially even further increase the advantage of iterative methods that use an initial guess from the previous point on an equilibrium curve; see also Section \ref{sec:Lyapunov}. Another important conclusion from the calculations in Figure \ref{fig:fig5} is that even though the two equilibrium curves have $\cO(10^{3})$ points each, the calculation took only a few seconds. Therefore the computation of all covariance matrices of equilibrium curves is expected to very fast on standard single-machine computer hardware for most small to medium-size ODE systems.

\begin{figure}[htbp]
\psfrag{a}{\scriptsize{(a)}}
\psfrag{b}{\scriptsize{(b)}}
\psfrag{b}{\scriptsize{(c)}}
\psfrag{iter}{\scriptsize{iter}}
\psfrag{fun}{\scriptsize{fun}}
\psfrag{time}{\scriptsize{time}}
\psfrag{tol}{\scriptsize{tol}}
	\centering
		\includegraphics[width=1\textwidth]{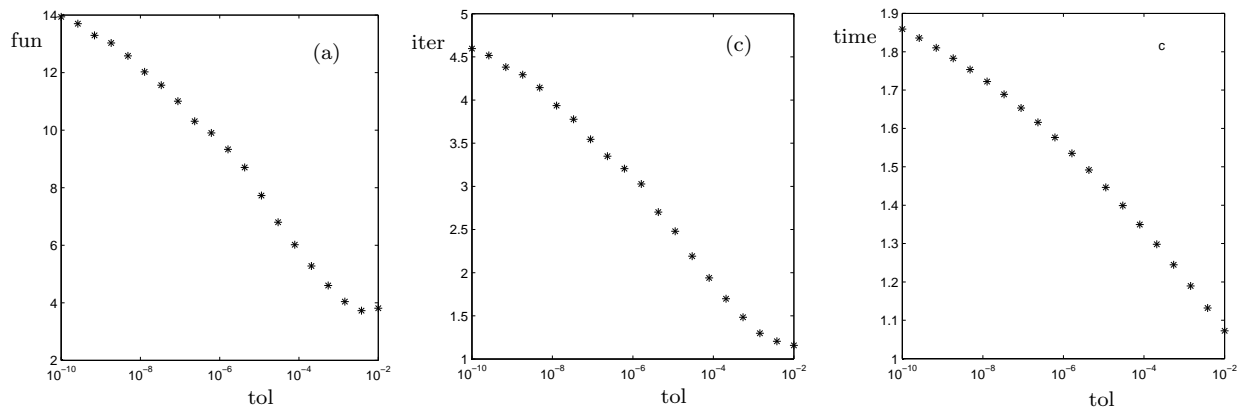}
	\caption{\label{fig:fig6}Computation of the distance between the two covariance ellipsoids of two stable equilibrium continuation curves (see Figure \ref{fig:fig4}) for given termination tolerances (logarithmic abscissa) of the SQP optimization algorithm. The distance is calculated at every tenth point equilibrium point; 216 distances have been computed in total. (a) Average number of function evaluations of the objective function over the 216 points required during the SQP algorithm. (b) Average number of iteration steps. (c) Total time in seconds to process 216 points.}
\end{figure} 

Figure \ref{fig:fig6} shows an overview of the computational cost to obtain the distances shown in Figure \ref{fig:fig4}(c) between ellipsoids using SQP as implemented in \cite{MatLab2010b}. The distance has been computed for 216 covariance ellipsoids sequentially along the equilibrium point curves. The initial conditions were obtained from the result of the previous optimization problem. The main result of Figure \ref{fig:fig6} is that the distance calculation can be carried out quickly and requires very few iterations steps and function evaluations. This means that we can evaluate the testfunction for overlapping ellipsoids efficiently using optimization. However, we do not claim that the algorithm we used here is optimal in any way. It is possible that other optimization techniques of methods to estimate distances between ellipsoids outperform the SQP approach we used here. However, from a practical point of view the results we obtain show that the computational time is certainly not prohibitive to process entire equilibrium bifurcation branches.\\

%%%%%%%%%%%%%%%%%%%%%%%%%%%%%%%%%%%%%%%%%%%%%%%%%%%%%
\section{A Predator-Prey System}
\label{sec:num_ex2}

In the previous section, we have focused on the algorithmic cost of our algorithm and the distance calculation between ellipsoids. In this section we are going to consider an example with a complicated noise term and focus on the value of our method for applications. The classical Rosenzweig-MacArthur \cite{RosenzweigMacArthur} model for the interaction of predators $Y$ and prey $X$ is given by
\be
\label{eq:RM}
\begin{array}{lcl}
x'&=&x\left(1-\frac{x}{\gamma}\right)-\frac{xy}{1+x}\\
y'&=&\beta\frac{xy}{1+x}-my
\end{array}
\ee
where $(x,y)$ represents the population densities of $(X,Y)$, $\gamma$ relates to the carrying capacity of the prey, $\beta$ is a conversion factor and $m$ a parameter describing mortality of the predator. The model \eqref{eq:RM} can be derived as a large-system size limit for the individual interactions between $X$ and $Y$. Finite-size effects of the population can be included into a stochastic fluctuation term. Using a Kramers-Moyal (or system-size) expansion one finds \cite{Roozen,Gardiner}
\be
\label{eq:RM_SDE}
\left(\begin{array}{c} dx_t\\ dy_t \\ \end{array}\right)= \left(\begin{array}{c}x\left(1-\frac{x}{\gamma}-\frac{xy}{1+x}\right) \\ \beta\frac{xy}{1+x}-my \end{array}\right) dt\\
+\sigma C(x,y)dW_t
\ee  
where $W_t=(W_t^{(1)},W_t^{(2)})^T$ is standard Brownian motion, the matrix-valued function $C$ is given by 
\benn
C(x,y)C(x,y)^T=B(x,y)\qquad \text{with}\quad
B(x,y)=\left(
\begin{array}{cc}
x\left(1+\frac{xy}{1+x}-\frac{x}{\gamma}\right) & -\frac{xy}{1+x} \\
-\frac{xy}{1+x} & y\left(\beta\frac{x}{1+x}+m\right)\\ 
\end{array}
\right)
\eenn
and $\sigma=\cO(1/N)$ where $N$ is the population size. Therefore $\sigma\ra 0$ corresponds to the limiting case of an infinite population which recovers the deterministic limit \eqref{eq:RM}. Observe that the noise terms in \eqref{eq:RM_SDE} are multiplicative and exhibit correlations between the two population densities. Therefore it is not immediately clear how a bifurcation diagram of \eqref{eq:RM} is altered once the (more realistic) finite-system size is considered.\\

\begin{figure}[htbp]
\psfrag{a}{(a)}
\psfrag{b}{(b)}
\psfrag{H}{H}
\psfrag{x}{$x$}
\psfrag{y}{$y$}
\psfrag{g}{$\gamma$}
	\centering
		\includegraphics[width=1\textwidth]{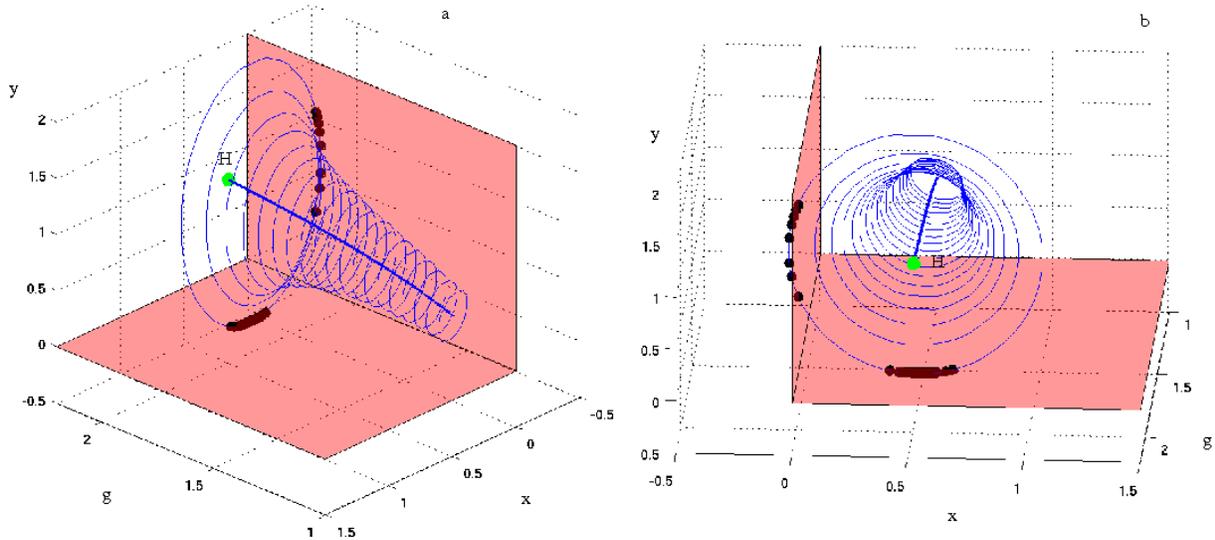}
	\caption{\label{fig:fig7}Continuation for \eqref{eq:RM_SDE} in $\gamma$ with $m=1$ and $\beta=3$ fixed. (a) and (b) are two viewpoints fo the same bifurcation diagram. The thick blue curve is computed via equilibrium continuation. The thin blue ellipses are computed with the algorithm from Section \ref{sec:main} for $h=1$ and $\sigma=0.01$. The Hopf bifurcation (H) at $\gamma=2$ is marked with a green dot. The red planes delimit the positive quadrant in $(x,y)$-phase space. The black dots on the ellipse at $\gamma=1.9$ are interpolation points for the ellipse outside the positive quadrant.}
\end{figure} 

We focus on deterministic Hopf bifurcations in the model which have received the most attention in ecological predator-prey models \cite{Kot}. Figure \ref{fig:fig7} shows an equilibrium continuation in $\gamma$ where increasing $\gamma$ can be interpreted as increasing the carrying capacity for the prey. Observe that a stable focus undergoes a Hopf bifurcation. Classical deterministic ecological theory \cite{Rosenzweig} argues that increasing the carrying capacity corresponds to enrichment and that the periodic solutions born in the Hopf bifurcation can move the system close to the coordinate axes
\benn
\{(x,y)\in\R^2:x\geq0\text{ and } y=0\} \qquad \text{and} \qquad \{(x,y)\in\R^2:y\geq0\text{ and } x=0\}
\eenn 
which delimit the positive quadrant. Once the system reaches any of the two axes it is easy to see that this corresponds to extinction of a species leading to a ``paradox of enrichment''. This ``paradox'' is a highly debated topic in ecology and many different ways of resolving it have been suggested, see for instance \cite{AbramsWalters,Jansen,Genkai-KatoYamamura,PetrovskiiLiMalchow}.

However, from our computation the ellipsoids suggest a very simple solution. The predator-prey system before a Hopf bifurcation can easily reach the axes as well, even for small noise which corresponds to a large (but finite!) population size. Close to the bifurcation point the ellipsoids increase in size which is precisely the well-known slowing down effect exploited in the theory of critical transitions \cite{KuehnCT1,KuehnCT2,Schefferetal}. If the carrying capacity in an ecosystem only increases slowly, which is reasonable to assume, then we expect that stochastic effects drive the system to extinction before the ``paradox of enrichment'' Hopf mechanism becomes relevant {i.e.} one would not see regular oscillations before extinction. Furthermore, the deterministic periodic solution occurring due to enrichment could actually have a stabilizing effect as the stochastic effects are small for a strongly attracting deterministic periodic orbit far from bifurcation. Indeed, the idea of stabilization of enrichment has been considered previously \cite{McCauleyMurdoch,Kirk}.   

%%%%%%%%%%%%%%%%%%%%%%%%%%%%%%%%%%%%%%%%%%%%%%%%%%%%%
\section{A Special Case - Kramer's Law}
\label{sec:Kramer}

So far, all computations only required local assumptions on the SDE \eqref{eq:SDE} regarding existence of a deterministic equilibrium and suitable smoothness. The covariance neigbhourhood $\cB(h)$ provides a rigorous local control of the dynamics. The global distance $\delta$ is a precisely computable, but probabilistically heuristic, measure to gain insight into global transition dynamics; without global assumptions on the dynamics this seems to be the best we can hope for. However, one may ask what happens if we have additional information on the global dynamics. Consider the SDE
\be
\label{eq:SDE_gradient}
dx_t=-\nabla U_\mu(x_t)dt+\sigma dW_t,\qquad \text{$x\in\R^n$, $U_\mu:\R^n\ra \R$, $W_t\in\R^n$, $\sigma\in\R$}
\ee
where the deterministic part is a gradient system with a potential $U_\mu$ parameterized by $\mu\in\R$. Critical points of $U_\mu$ correspond to equilibria for the deterministic dynamics. Fix some $\mu\in\R$ and suppose $U_\mu$ has precisely two local minima $x^{*}$ and $y^*$, corresponding to stable equilibria, and one saddle point $z^*$. Define
\benn
\tau_{\cN(y^*)}:=\inf\{t>0:x_t\in\cN(y^*)\},\qquad x_0=x^*, 
\eenn 
for a suitable neighbourhood $\cN(y^*)$ of $y^*$. Under the assumption that the saddle point $z^*$ has a single unstable eigendirection with eigenvalue $\lambda(z^*;\mu)>0$, the Eyring-Kramers law \cite{Eyring,Kramers} states that
\be
\label{eq:Kramers_addref}
\E[\tau_{\cN_{\sigma^2/2}(y^*)}|x_0=x^*]=\frac{2\pi}{|\lambda(z^*;\mu)|}\sqrt{\frac{|\det(A(z^*;\mu))|}{\det(A(x^*;\mu))}}e^{2[U_\mu(z^*)-U_\mu(x^*)]/\sigma^2}\left(1+\cO(\sigma |\ln(\sigma^2/2)|^{3/2})\right)
\ee
where $A(x^*;\mu)=D^2U_\mu\in\R^{n\times n}$ is the Hessian of $U_\mu$ and $\cN_{\sigma^2/2}(y^*)$ a ball of radius $\sigma^2/2$ around $y^*$. The precise formula \eqref{eq:Kramers_addref} is due to Bovier et al. \cite{BovierEckhoffGayrardKlein}; see also \cite{HaenggiTalknerBorkovec,Berglund3} for reviews and generalizations of Kramers' law. The probability of switching due to noise from $x^*$ to $y^*$ is given to leading-order by \eqref{eq:Kramers_addref} and interchanging the roles of $x^*$ and $y^*$ provides the noise-induced switching estimates for the transition from $y^*$ to $x^*$. To compute \eqref{eq:Kramers_addref} we can follow an analogous strategy as for the more general case discussed so far. The equilibria $x^*$, $y^*$ and $z^*$ as well as the associated linearizations $A(\cdot;\mu)$ can be computed efficiently via numerical continuation for a curve parameterized by $\mu\in\R$, the function $V$ is available by assumption and it remains to compute $\lambda(z^*;\mu)$, $\det(A(z^*;\mu))$ and $\det(A(x^*;\mu))$. To compute the determinants we can simply use the LU decomposition \cite{GolubVanLoan} which also works well for large sparse systems. However, calculating the leading eigenvalue is bound to be costly if we look to compute all eigenvalues and then extract the leading one. Suppose we are given the results $A(z^*;\mu_{k-1})$, $\lambda(z^*;\mu_{k-1})$ and the associated eigenvector $v(z^*;\mu_{k-1})$ from at the last continuation step then we can again use an iterative method to compute $\lambda(z^*;\mu_{k})$. For example, setting $v_0(z^*;\mu_{k})=v(z^*;\mu_{k-1})$ and $\lambda_0(z^*;\mu_{k})=\lambda(z^*;\mu_{k-1})$ then Rayleigh quotient iteration \cite{Parlett} is given by
\be
\label{eq:Rayleigh1}
v_{j+1}(z^*;\mu_{k})=\frac{(A(z^*;\mu_{k})-\lambda_{j}(z^*;\mu_{k})I)^{-1}v_j(z^*;\mu_{k})}{\|(A(z^*;\mu_{k})-\lambda_{j}(z^*;\mu_{k})I)^{-1}v_j(z^*;\mu_{k})\|}, \qquad \text{for $j=0,1,2,\ldots$}
\ee     
and the eigenvalue for the $j$-th iteration step is 
\be
\label{eq:Rayleigh2}
\lambda_j(z^*;\mu_k)=\frac{v_j(z^*;\mu_k)^TA(z^*;\mu_k)v_j(z^*;\mu_k)}{v_j(z^*;\mu_k)^Tv_j(z^*;\mu_k)}.
\ee
It is well-known that for a symmetric matrix $A(z^*;\mu_j)$ the iteration \eqref{eq:Rayleigh1}-\eqref{eq:Rayleigh2} converges cubically to the leading eigenvalue and eigenvector \cite{Parlett,BattersonSmillie1}; in particular, $\lambda_j(z^*;\mu_k)\ra \lambda(z^*;\mu_k)$. Since $A(z^*;\mu_j)=D^2U_{\mu_j}$ is derived from a, sufficiently smooth, potential $U_{\mu_j}$ the matrix at each continuation step is symmetric and the fast convergence results for Rayleigh iteration apply; note that this may not be the case for open sets of ``bad'' starting conditions if the matrix is not symmetric \cite{BattersonSmillie}. In any case, evaluating the remaining terms in \eqref{eq:Kramers_addref} is straighforward so that we can calculate mean-first passage times between equilibria in gradient systems quickly, with high accuracy, and rigorous error estimates by using numerical continuation.\\

\textbf{Acknowledgments:} I would like to thank Tilo Schwalger for insightful discussions about noise in neuronal models and John Guckenheimer, Daniele Avitabile and Thorsten Riess for interesting discussions about an earlier draft of this paper. Furthermore, comments of two anonymous referees helped to improve the manuscript.

\bibliographystyle{plain}
\bibliography{../my_refs}

\end{document}